\documentclass[10pt,reqno]{amsart}

\usepackage{amsmath,amscd}

\usepackage[sort]{natbib}

 \newcommand\CC{\mathcal C}

 \newcommand\D{\mathcal D}

\newcommand\G{ G} \renewcommand\H{H} 
 \newcommand\PP{\mathbb P}
 
\newcommand\Y{\mathcal Y} 
 \newcommand\OO{\mathcal O}

\newcommand{\mult}{\operatorname{mult}}

\makeatletter \@addtoreset{equation}{section} \makeatother

\newtheorem{thm}[equation]{Theorem} \newtheorem*{thm*}{Theorem}
\newtheorem{lem}[equation]{Lemma}
\newtheorem{cor}[equation]{Corollary}
\newtheorem{prop}[equation]{Proposition}
\newtheorem{conj}[equation]{Conjecture} \theoremstyle{definition}

\newtheorem{rmk}[equation]{Remark} \newtheorem*{rmk*}{Remark}

\title[Seshadri constants on  symmetric products of curves]{Seshadri constants on\\
  symmetric products of curves}
\author{J. Ross}

\begin{document}
\bibliographystyle{abbrv}


\begin{abstract}
  Let $X_g=C^{(2)}_g$ be the second symmetric product of a very
  general curve of genus $g$.  We reduce the problem of describing the
  ample cone on $X_g$ to a problem involving the Seshadri constant of a
  point on $X_{g-1}$.  Using this we recover a result of
  Ciliberto-Kouvidakis that reduces finding the ample cone of $X_g$
  to the Nagata conjecture when $g\ge 9$. We also give new bounds on
  the the ample cone of $X_g$ when $g=5$.
\end{abstract}

\maketitle
                              
\section{Introduction}

Consider the second symmetric product $C^{(2)}$ of a smooth curve $C$
of genus $g\ge 2$.  This smooth surface comes with some naturally defined
divisors.  Given a point $p\in C$ there is the divisor $x_p=\{p+q :
q\in C\}$ whose numerical class is independent of $p$ and will be
denoted by $x$.  Another divisor on $C^{(2)}$ is given by the diagonal
$\Delta=\{p+p | p\in C\}$ whose numerical class is denoted by
$\delta$.


We will be interested in describing the the intersection $N$ of the
ample cone with the plane in $N^1(C^{(2)})_{\mathbb R}$ spanned by $x$
and $\delta$.  Note that when $C$ is a very general curve the classes
$x$ and $\delta/2$ generate $N^1(C^{(2)})$ so in this case $N$ is the
entire ample cone of $C^{(2)}$.  Since $N$ is a two dimensional cone it
is described by two boundary rays.  The first boundary is easily
given: since the diagonal $\Delta$ is an irreducible curve of negative
self-intersection it spans a boundary of the effective cone,
so its dual ray is one boundary.  The more interesting boundary of $N$
is characterised by the quantity
\begin{eqnarray*}
  \tau(C) &=& \text{inf}\,\, \{ s>0 \,: (s+1)x-(\delta/2) \text{ is ample}\,\}.
\end{eqnarray*}
There is the obvious universal bound
$$\tau(C)\ge \sqrt{g},$$
coming from the fact that if $(s+1)x-(\delta/2)$ is ample then it has
positive self-intersection. The following conjecture governs the nef
cone of a $C^{(2)}$ when $C$ is very general.

\begin{conj}\label{conj:ampleconeoncurve}
  If $C$ is a very general curve of genus $g\ge 4$ then
  $\tau(C)=\sqrt{g}$.
\end{conj}

This conjecture asserts that for a very general curve $C$ the other boundary of $N$ has zero-self intersection, and thus the nef (resp.\ ample) cone of $C^{(2)}$ is as large as possible.  It is only known to hold when $g$ is a perfect square  \cite{kouvidakis:93:divis_symmet_produc_curves}.    \\

Our aim is to give a lower bound for $\tau(C)$ in terms of the
Seshadri constant of a point in $D^{(2)}$ where $D$ is a smooth curve of
genus $g-1$.  If $X$ is a smooth surface and $L$ is the (numerical
class of) a nef $\mathbb R$-divisor on $X$ the Seshadri constant at a
collection of distinct points $p_1,\ldots, p_m\in X$ is defined to be
$$ \epsilon(p_1,\ldots,p_m;X,L) = \text{inf}_C\left\{ \frac{L.C}{\sum_i \mult_{p_i} C}\right\},$$
where the infimum is over all reduced irreducible curves $C\subset X$ passing
through at least one of the $p_i$.  We will prove the following
connecting Seshadri constants and the ample cone of second
symmetric products of smooth curves.

\begin{thm}\label{thm:ampleandseshadristronger}
  Let $D$ be a smooth curve of genus $g-1$.  Suppose $a,b>0$ are such
  that $a/b>\tau(D)$ and for a very general point $p\in D^{(2)}$
$$\epsilon\left(p; D^{(2)}, (a+b )x - b (\delta/2)\right)\ge b.$$
Then for a very general curve $C$ of genus $g$,
$$\tau(C) \le \frac{a}{b}.$$
\end{thm}
Thus Conjecture
\ref{conj:ampleconeoncurve} is implied by the following conjecture
about Seshadri constants:
\begin{conj}\label{conj:seshadriconstantismaximal}
  Let $D$ be a very general curve of genus $g-1$ with $g\ge 5$ and $p$
  be a very general point in $D^{(2)}$.  Then
  \begin{equation}
\epsilon\left(p; D^{(2)},(\sqrt{g} +1 )x - (\delta/2)\right)=1. \label{eq:seshadrequality}
\end{equation}
\end{conj}

We note that this conjecture is not easy to prove as the class
$L=(\sqrt{g}+1)x-(\delta/2)$ has degree $L^2=1$ so the equality in
\eqref{eq:seshadrequality} asserts that the Seshadri constant of this
$\mathbb R$-divisor is maximal (see \ref{eq:seshadriupperbound}).
There is a general lower bound due to Ein-Lazarsfeld
\cite{ein-lazarsfeld:93:seshad_const_smoot_surfac} for the Seshadri
constant of general points in surfaces with respect to
\emph{integral} divisors but this does not extend to the case of
$\mathbb R$-divisors.  However when $g$ is a perfect square we can
apply \cite{ein-lazarsfeld:93:seshad_const_smoot_surfac} to deduce
that $\epsilon(p;D^{(2)},L)= 1$ for a very general $p\in D^{(2)}$
where $D$ is a smooth curve of genus $g-1$.  Thus we get another proof
that if $g\ge 4$ is a perfect square then $\tau(C) =\sqrt{g}$ for a
very general curve $C$ of genus $g$.

  \begin{rmk}
    As pointed out by Lazarsfeld, it is not the case that the analogy
    of Conjecture \ref{conj:seshadriconstantismaximal} holds for all
    polarised surfaces $(X,L)$ such that $L$ is an ample $\mathbb
    R$-divisor with $L^2=1$ and $L.E\ge 1$ for all but finitely many
    curves $E\subset X$.  For example let $X$ be an abelian surface of
    type $(1,d)$ with Picard number 1 generated by the ample line
    bundle $L'$ and set $L=L'/\sqrt{2d}$.  Then $L^2=1$ and $L.E\ge 1$
    for all irreducible curves $E\subset X$.  But whenever $\sqrt{2d}$
    is irrational it is known that $\epsilon(p,X;L')<\sqrt{2d}$ (in
    fact it is rational \cite{bauer(99):seshad_const_algeb_surfac}) so
    $\epsilon(p;X,L)<1$ for all $p\in X$.
  \end{rmk}

  The proof of Theorem \ref{thm:ampleandseshadristronger} uses a
  degeneration of the symmetric product that arises from letting $C$
  degenerate to the nodal curve $C_0$ obtained by gluing two points
  in $D$.  The same degeneration allows us to compare the multipoint
  Seshadri constants of $C^{(2)}$ and $D^{(2)}$.  We define
  $\epsilon_m(X,L)$ to be the Seshadri constant of a collection of $m$
  very general points in $X$.

\begin{thm}\label{thm:seshadrimulti}
  Let $D$ be a smooth curve of genus $g-1$ and fix an integer $m\ge
  1$.  Suppose that there are numbers $a,b>0$ with $a/b>\tau(D)$ such
  that
$$\epsilon_{m+1}\left( D^{(2)}, (a+b)x-b(\delta/2)\right) \ge b.$$
Then for a very general curve $C$ of genus $g$ the class $
(a+b)x-b(\delta/2)\in N^1(C^{(2)})$ is nef and
$$\epsilon_{m} \left(C^{(2)},  (a+b)x-b(\delta/2)\right) \ge b.$$
\end{thm}

A more concise way to state this theorem is to let $\epsilon_{m,g}(s)$
be the Seshadri constant of $m$ very general points in $C^{(2)}$ with
respect to the class $(s+1)x-\delta/2$, where $C$ is a very general
curve of genus $g\ge 0$.  Then
$$\epsilon_{m,g}(s) \ge \epsilon_{m+1,g-1}(s),$$
where this is to be interpreted as holding for all $s,g$ such that the
right hand makes sense.  \\

As the second symmetric product of $\PP^1$ is $\PP^2$, induction on
$g$ yields:

\begin{cor}[Ciliberto-Kouvidakis \cite{ciliberto-kouvidakis(99):symmet_produc_a_curve_with_gener_modul}]\label{cor:ampleintermsofseshadriofpp2}
Let $C$ be a very general curve of genus $g\ge 1$.  Then 
$$\tau(C) \le \frac{1}{\epsilon_g(\PP^2,\OO_{\PP^2}(1))}.$$
\end{cor}

As is well known, one formulation of the Nagata conjecture states that
if $g\ge 9$ then the Seshadri constant of $g\ge 9$ very general points
in $\PP^2$ is maximal:

\begin{conj}[Nagata Conjecture]\label{conj:nagata}
  If $g\ge 9$ then
$$\epsilon_g(\PP^2,\OO_{\PP^2}(1)) = \frac{1}{\sqrt{g}}.$$
\end{conj}
Thus, as proved in
\cite{ciliberto-kouvidakis(99):symmet_produc_a_curve_with_gener_modul},
the Nagata conjecture yields the ample cone of $C^{(2)}$ for a very
general curve of genus $g\ge 9$.  Currently the Nagata conjecture is
only proved when $g$ is a perfect square.  For other $g$ there are
several bounds on $\epsilon_g(\PP^2,\OO_{\PP^2(1)})$ (e.g.
\cite{szemberg_tutaj_gasinska(02):gener_blow_ups_projec_plane,xu(94):curves_in_f_p_sp,strycharz_szemberg_szemberg(04):remar_nagat_conjec,harbourne_roe:05:multip_seshad_const_bbb_p_sp_2}).
For instance using \eqref{cor:ampleintermsofseshadriofpp2} and results
from \cite{strycharz_szemberg_szemberg(04):remar_nagat_conjec} we get
that for a very general curve of genus $g\ge 10$,
\begin{equation}
\tau(C) \le \frac{\sqrt{g}}{\sqrt{1- \frac{1}{g+1}}}.\label{eq:boundfromwhatisknownaboutnagata}
\end{equation}
\noindent When $g$ is not a perfect square this improves on the bound
$\tau(C)\le \frac{g}{[\sqrt{g}]}$ from
\cite{kouvidakis:93:divis_symmet_produc_curves}. 

\begin{rmk}\label{rmk:conjecturerelations}
From the discussion above there are implications \\

\hspace{.6cm}\begin{minipage}{3cm}
Conjecture \ref{conj:nagata} \\ 
{\small \it{(Nagata conjecture)}}
\end{minipage} 
$\Rightarrow\quad$
\begin{minipage}{3cm}
\text{Conjecture } \ref{conj:seshadriconstantismaximal}\\
{\small \it{(Seshadri constants)}}
\end{minipage} 
$\Rightarrow\quad$
\begin{minipage}{3cm}
\text{Conjecture } \ref{conj:ampleconeoncurve}\\
{\small {\it (Ample cone of $C^{(2)}$)}}
\end{minipage} \\

\noindent so Conjecture \ref{conj:seshadriconstantismaximal}
concerning Seshadri constants sits between the Nagata conjecture and
the conjecture governing the ample cone of a general $C^{(2)}$.  It is
possible that Conjecture \ref{conj:seshadriconstantismaximal} is
easier than the full Nagata conjecture but there is, of course,
currently no proof of this.
\end{rmk}

\begin{rmk}
  The degeneration we use in the proofs of the above theorems is
  related to the Franchetta degeneration used in
  \cite{ciliberto-kouvidakis(99):symmet_produc_a_curve_with_gener_modul}
  which arises from degenerating $C$ to a rational nodal curve.  The
  main difference is that here we consider the case that $C$ develops
  one node at a time.  Moreover rather than using a degeneration of
  $C^{(2)}$ we find it easier to use a degeneration of $C\times C$ and
  only consider divisors that are invariant under permuting the
  factors.  A related degeneration of of $C^{(2)}$ coming from letting
  $C$ develop cusps is described in
  \cite{pacienza_polizzi:01:a_degen_symmet_produc_a}.
\end{rmk}

The values of $\tau(C)$ are known for a very general curve of genus
$g\le 4$ (see Section \ref{sec:lowgenus}) and from the discussion
above the Nagata conjecture governs the case that $g\ge 9$.  There
appears to be little known in the intermediate range $5\le g\le 8$ and
one would imagine that one of the two extreme cases holds (namely that
either $\tau(C) =\epsilon_g(\PP^2,\OO_{\PP^2}(1))^{-1}$ or
$\tau(C)=\sqrt{g}$ for a very general curve $C$ of genus $5\le g\le
8$).  In Section \ref{sec:applicationgenus5} we apply Theorem
\ref{thm:ampleandseshadristronger} to show that the first case does
not hold when $g=5$; more precisely we show that if $C$ is a very
general curve of genus $5$ then $\tau(C)\le 16/7$ which gives
$$2.236\simeq \sqrt{5}\le \tau(C)\le 16/7\simeq 2.286< \epsilon_5(\PP^2,\OO_{\PP^2})^{-1} = 2.5.$$
To achieve this we use the techniques of Ein-Lazarsfeld
\cite{ein-lazarsfeld:93:seshad_const_smoot_surfac} to get lower bounds
of Seshadri constants of points in $D^{(2)}$ where $D$ is a curve of
genus $4$.  This bound is stronger than that obtained from Corollary
\ref{cor:ampleintermsofseshadriofpp2} as
$\epsilon_5(\PP^2,\OO_{\PP^2}(1))=2/5$
\cite{strycharz_szemberg_szemberg(04):remar_nagat_conjec}.  The number
$16/7$ is not expected to be optimal, but it is, as far as I am aware,
the best that is currently known.\vspace{4mm}

\noindent {\bf Acknowledgements:} I would like to thank Robert
Lazarsfeld for helpful discussions and for pointing out that the hope
that the analogy of Conjecture \ref{conj:seshadriconstantismaximal}
holds for all surfaces was too optimistic.  I also thank Ciro
Ciliberto, Johan de Jong, Ian Morrison and Yusuf Mustopa.\\

\noindent {\bf Notation and conventions:} We work throughout over
$\mathbb C$.  The N\'eron-Severi space of divisors (resp.\ $\mathbb
R$-divisors) on a variety $V$ modulo numerical equivalence is denoted
$N^1(V)$ (resp.\ $N^1(V)_{\mathbb R}$).  An $\mathbb R$-divisor $L$ on
a variety $V$ is ample (resp.\ nef) if it is a formal sum
$\sum_{i=1}^r a_i D_i$ of ample (resp.\ nef) divisors where the $a_i\in \mathbb R$
are positive (resp.\ non-negative).  Equivalently $D\in
N^1(V)_{\mathbb R}$ is nef if and only if it has non-negative self
intersection with every irreducible curve $C\subset V$.

We say that $p\in V$ is a very general point if there is a countable
collection of proper subvarieties $(V_n)_{n\ge 1}$ of $V$ such that
$p$ is not contained in the union $\bigcup_{n\ge 1} V_n$.  A
collection $p_1,\ldots,p_m$ of points in $V$ is very general if
$(p_1,\ldots,p_m)\in V^{\times m}$ is very general.  By a very general
curve we mean a smooth curve whose corresponding point in the moduli
space $M_g$ is very general.


\section{Preliminaries}\label{sec:preliminaries}
%

\subsection{Divisors the second symmetric product.}


Let $C$ be a smooth curve of genus $g\ge 0$.  The product $C\times C$
has a natural involution and the second symmetric product is the
quotient $\sigma_C\colon C\times C\to C^{(2)}$ which is a smooth
surface.  We denote the image of a point $(p,q)\in C\times C$ by
$p+q$.  In $N^1(C^{(2)})$ we have the classes $x$ and $\delta$ as
defined in the introduction. (Of course these classes really depend on
$C$ but this will always be clear from context.) It is well known that
when $C$ is a very general curve then $N^1(C^{(2)})$ is spanned by $x$
and $\delta/2$
(\cite{arbarello_cornalba_griffiths(85):geomet_algeb_curves} p.359)
and when $g\ge 1$ they are independent.  When $g=0$, $C^{(2)}=\PP^2$
and both $x$ and $(\delta/2)$ are the class of the hyperplane.

The intersection of these classes is given by $x^2=1$, $\delta^2=4-4g$,
and $x.\delta=2$, so
$$
\left((n+\gamma) x - \gamma (\delta/2)\right) \cdot \left((n'+\gamma') x - \gamma'
(\delta/2)\right) = nn'-\gamma\gamma' g.$$ Notice that if $(n+\gamma)x -
\gamma(\delta/2)$ is effective then intersecting with the ample class
$x$ implies $n>0$.

As in the introduction define
\begin{eqnarray*}
  \tau(C) &=& \inf\,\, \{ s>0 \,: (s+1)x-(\delta/2) \text{ is ample}\,\}\\  &=&\text{min} \{ s\ge 0 : (s+1)x-(\delta/2) \text{ is nef}\,\}.
\end{eqnarray*}
Since $(\tau(C)+1)x - (\delta/2)$ is nef, it has non-negative
self-intersection which yields $$\tau(C)\ge \sqrt{g}.$$ Now as is
standard in such situations, the function $\tau(C)$ is
semicontinuous with respect to $C$:

\begin{lem}\label{lem:jumpsdown}
  Let $X\to T$ be a flat family of smooth curves over an irreducible
  base $T$ and for $t\in T$ denote the fibre by $C_t$.  If $t_0\in T$ is
  fixed then
  $$\tau(C_{t})\le \tau(C_{t_0}) \text{ for very general } t\in T.$$
\end{lem}
\begin{proof}
  Let $Y\to T$ be the relative second symmetric product of $X$, and
  denote the fibre of $Y$ over $t$ by $Y_t$.  Let $\D\subset Y$ be the
  diagonal and pick a divisor $D$ on $Y$ which has class $x$ on a each
  fibre $Y_t$.  If $\tau=\tau(C_{t_0})$ then by hypothesis
  $F=(\tau+1)D-(\mathcal D/2)$ restricts to a nef divisor on $Y_{t_0}$.
  Hence for very general $t$ the restriction if $F$ to $Y_{t}$ is
  nef (\cite{lazarsfeld(04):posit_in_algeb_geomet} 1.4.14) which
  implies $\tau(C_{t})\le \tau$.
\end{proof}

In particular by applying this to a complete family of smooth curves
we see that if $\tau(C)\le \tau_0$ for some smooth curve $C$ of genus
$g$, then the same bound holds for a very general curve of genus $g$.
Moreover if $C$ is a very general curve then $\tau(C)$ is independent
of the actual curve chosen. \\

A geometric interpretation of Conjecture \ref{conj:ampleconeoncurve}
can be given in terms of the existence of ``exceptional'' curves in
$C^{(2)}$:

\begin{lem} Let $C$ be a smooth curve of genus $g\ge 2$.\label{lem:exceptionalcurvesonC2}
  \begin{enumerate}
  \item If $\tau(C)>\sqrt{g}$ then there exists a reduced irreducible
    curve $D\subset C^{(2)}$ with numerical class
    $(n+\gamma)x-\gamma(\delta/2) + \sigma$ where $\sigma.x=\sigma.\delta=0$
    such that $\tau(C) = \frac{\gamma g}{n}$.

  \item If $C$ is a very general curve then $\tau(C)=\sqrt{g}$ if and
    only if $\Delta$ is the only reduced irreducible curve in $C^{(2)}$ with
    negative self-intersection.  Moreover if there does exist another such
    curve of negative self intersection then it is unique.
  \end{enumerate}
\end{lem}

\begin{proof}
   Suppose $\tau=\tau(C)>\sqrt{g}$.  Then the $\mathbb R$-divisor
  $$F=(\tau+1)x-(\delta/2)$$
  has positive self-intersection and, by definition of $\tau$, is nef
  but not ample. Thus by the Nakai criterion for real divisors
  \cite{campana_peternell(90):algeb_ample_cone_projec_variet} there is
  a reduced irreducible curve $D\subset C^{(2)}$ with $D.F=0$.  We can write
  the numerical class of $D$ as $(n+\gamma)x-\gamma (\delta/2)+\sigma$ where
  $\sigma$ is a class orthogonal to $x$ and $\delta$ so $D.F=0$
  implies $\tau=\frac{\gamma g}{n}$.  We note that since  $n>0$
  this implies $\gamma >0$.

  Now suppose $C$ is very general.  Then the effective cone of
  $C^{(2)}$ is spanned by $x$ and $\delta/2$. Thus if
  $\tau=\tau(C)>\sqrt{g}$ the curve $D$ above has class
  $(n+\gamma)x-\gamma(\delta/2)$ (i.e.\ $\sigma=0$).  Hence
  $D^2=n^2-g\gamma ^2 = \frac{n^2}{g}(g-\tau^2)<0$ so $D$ has negative
  self intersection and clearly $D\neq \Delta$ as $\gamma >0$.

  Before proving the converse we deal with uniqueness.  To this end
  suppose that $(n+\gamma)x - \gamma(\delta/2)$ and $(n'+\gamma')x -
  \gamma'(\delta/2)$ are classes of distinct reduced irreducible
  curves of negative self intersection.  Then $n^2-\gamma^2 g<0$,
  $n'^2-\gamma'^2<0$ and $nn'-\gamma\gamma'g\ge 0$ which implies
  $\gamma_1$ and $\gamma_2$ have opposite sign.  Hence any irreducible
  curve $D\neq \Delta$ in $C^{(2)}$ with $D^2<0$ must have numerical
  class $(n+\gamma)x - \gamma(\delta/2)$ with $\gamma>0$, and if it
  exists it is unique.  Thus if there exists a reduced irreducible
  curve $D\neq \Delta$ with negative self intersection it has numerical
  class $(n+\gamma)x-\gamma(\delta/2)$ with $\gamma>0$ and
  $D^2=n^2-g\gamma^2<0$.  As $F$ is nef we know $0\le F.D = \tau n -g
  \gamma$ which implies $\tau\ge \frac{g\gamma}{n}>\sqrt{g}$.
\end{proof}

Before proceeding with the main results of this paper we digress to
discuss a finiteness result concerning the possible values of
$\tau(C)$ as $C$ ranges over all curves of fixed genus $g\ge 2$.

\begin{prop}\label{prop:finitenessofpossibletau}
  Fix a real number $\alpha>\sqrt{g}$.  Then 
$$\{\tau\ge \alpha: \tau=\tau(C) \text{ for some smooth curve } C \text{ of genus  }g \}$$
is a finite set.  Equivalently the only possible accumulation point of
the set given by  $\{\tau(C) : C \text{ a smooth curve of genus }g\}$ is $\sqrt{g}$.
\end{prop}
\begin{proof}
  Fix a number $s\in (\sqrt{g},\alpha)\cap \mathbb Q$.  We first prove
  that there exists an integer $k$ such that for any smooth curve $C$ of
  genus $g$ the divisor $k[(s+1)x-(\delta/2)]$ on $C^{(2)}$ is
  effective.

  To this end let $C$ be any smooth curve of genus $g$ and fix a
  $\mathbb Q$-divisor $F$ on $C^{(2)}$ whose numerical class is
  $(s+1)x-(\delta/2)$.  We note that the canonical class of $C^{(2)}$
  is $K=(2g-2)x-(\delta/2)$
  (\cite{kouvidakis:02:some_resul_morit_their_applic} Prop. 2.6).
  Thus if $k\in\mathbb N$ is sufficiently large (with $ks\in \mathbb
  N$) then $x.(K-kF)<0$. So by Serre duality $h^2(\OO(kF)) =
  h^0(\OO(K-kF)) =0$ as $x$ is ample.  Hence for such $k$,
$$h^0(\OO(kF)) \ge h^0(\OO(kF)) -h^1(\OO(kF)) = \chi(kF) = p(k)$$
where by the Riemann-Roch theorem $p(k)$ is a polynomial whose
coefficients depend only on $s$ and $g$ (and not on the specific curve
$C$ or choice of $F$).  The leading order coefficient of $p(k)$ is $F^2/2 =
(s^2-g)/2>0$ so there exists a $k$ (independent of $C$) such
that $h^0(\OO(kF))>0$ as claimed.

Now suppose $C$ is chosen so that $\tau(C) \ge \alpha$.  By Lemma
\ref{lem:exceptionalcurvesonC2}(a) there exists a reduced irreducible
curve $D\subset C^{(2)}$ with numerical class $(n+\gamma)x -
\gamma(\delta/2)+ \sigma$ where $\sigma.x=\sigma.\delta=0$ and
$\tau(C)= \frac{g\gamma}{n}$.  With $k,F$ as above there is a divisor
$E\subset |kF|$.  But as $s<\alpha\le \tau(C)$ we have $D.F =
(ns-\gamma g) <0$ which implies that $D\subset E$.  Since $C$ is reduced
and irreducible this implies that $n=x.D\le x.E = ks$.  Thus
letting $N:=ks$ we get that $n\le N$.

To complete the proof we use that the divisor $G=(g-1)x+(\delta/2)$ is
always nef (as it is dual to the diagonal).  Hence $D\subset E$ also
implies $D.G\le E.G$ which yields $gn+\gamma g \le k(gs+g)$ so $\gamma
\le k(s+1)=:M$.  Thus $\tau(C)$ lies in the set
$$\{ \tau : \tau=\frac{g\gamma}{n} \text{ with } n,\gamma\in \mathbb N\ \text{ and } n\le N, \gamma \le M\}$$
which is finite.

\end{proof}

\begin{rmk}
  A similar finiteness result for Seshadri constants in families of
  surfaces can be found in
  \cite{oguiso(02):seshad_const_in_a_famil_surfac}.  It would be
  interesting to know if the finiteness from Proposition
  \ref{prop:finitenessofpossibletau} still holds when
  $\alpha=\sqrt{g}$.

\end{rmk}

\subsection{The case of low genus}\label{sec:lowgenus}

For low genus it is possible to describe the intersection of the ample
cone with the plane spanned by $x$ and $\delta$ by finding explicit
irreducible curves of negative self intersection.  For details see
\cite{ciliberto-kouvidakis(99):symmet_produc_a_curve_with_gener_modul,kouvidakis:93:divis_symmet_produc_curves}
(or \cite{lazarsfeld(04):posit_in_algeb_geomet} Section 1.5.B).

\begin{itemize}
\item $g=0$: Here $C^{(2)}=\PP^2$ and $(s+1)x-(\delta/2)=sh$ where $h$
  is the class of the hyperplane, so trivially $\tau(\PP^1) = 0$.
\item $g=1$: In this case it is well known that if $C$ is a very
  general genus 1 curve then the closure of the effective cone of
  $C^{(2)}$ is the nef cone.  It is a closed circular cone described
  by the equations $\alpha^2\ge 0$, $\alpha.h\ge 0$ where $h$ is an
  ample class (\cite{lazarsfeld(04):posit_in_algeb_geomet} Lemma
  1.5.4).  Thus $\tau(C)=1$.
\item $g=2$: Any curve $C$ of genus $2$ is hyperelliptic.  Using the
  $g_1^2$ one can produce an irreducible curve in $C^{(2)}$ of negative self
  intersection whose class is $2x-(\delta/2)$, and thus $\tau(C)=2$.
\item $g=3$: If $C$ is a very general curve of genus $2$ then it is
  possible to construct an irreducible curve in $C^{(2)}$ whose class
  is $16x-6(\delta/2)$ and thus has self-intersection -8
  \cite{kouvidakis:93:divis_symmet_produc_curves,ciliberto-kouvidakis(99):symmet_produc_a_curve_with_gener_modul}.
  Using this one deduces that $\tau(C) = 9/5$.
\item $g=4$: If $C$ is a very general curve of genus $4$ then
  $\tau(C)=2$.  In fact any such curve admits two $g_1^3$, and the
  associated $\Gamma_3$ is an irreducible curve whose class is
  $3x-(\delta/2)$ spans one boundary of the effective cone.  (This can
  also be obtained from Corollary
  \ref{cor:ampleintermsofseshadriofpp2}).
\end{itemize}

\subsection{Seshadri constants}\label{sec:seshadri}

We record some basic definitions and properties of Seshadri constant
and refer the reader to
\cite{bauer(99):seshad_const_algeb_surfac,lazarsfeld(04):posit_in_algeb_geomet}
for a comprehensive treatment.  Let $X$ be a smooth variety of
dimension $n$ and $L$ be a nef numerical class in $N^1(X)_{\mathbb
  R}$.  If $p_1,\ldots,p_m$ are points in $X$ define
$$\epsilon(p_1,\ldots,p_m; X,L) = \text{inf}_C \left\{ \frac{L.C}{\sum_{i=1}^r \mult_{p_i} C}\right\},$$
where the infimum is over all reduced irreducible curves $C$ in $X$ that pass
through at least one of the $p_i$. Equivalently if $\pi\colon B\to X$
is the blowup of $X$ at these points with exceptional divisor $E$
then
$$\epsilon(p_1,\ldots,p_m;X,L) = \text{max}\{s\ge 0 : \pi^* L - sE \text{ is nef}\}.$$
By a standard semicontinuity argument similar to \eqref{lem:jumpsdown}
the Seshadri constant of $m$ very general points does not depend on
the actual points chosen.  Thus we can set
$$\epsilon_m(X,L) = \epsilon(p_1,\ldots,p_m;X,L) \text{ where } p_1,\ldots,p_m \text{ are in very general position}.$$
Notice that if $\pi^* L - cE$ is nef then $(\pi^* L - cE)^n\ge 0$
which implies
\begin{equation}
 \epsilon(p_1,\ldots,p_m;X,L)\le \sqrt[n]{\frac{L^n}{m}}.\label{eq:seshadriupperbound} 
\end{equation}
It is an extremely interesting and difficult problem to get lower
bounds for Seshadri constants in general
\cite{ein-lazarsfeld:93:seshad_const_smoot_surfac,ein_kuechle_lazarsfeld(95):local_posit_ample_line_bundl}
or even to calculate them in examples (see
\cite{lazarsfeld(04):posit_in_algeb_geomet} and the reference
therein).

We will make use of the following simple lemma which says that Seshadri
constants of a collection of points in $C^{(2)}$ can be calculated by
looking at their preimage in $C^{\times 2}$.

\begin{lem}\label{lem:seshadriquotient}
  Let $D$ be a smooth curve and $\sigma_D \colon D^{\times 2} \to
  D^{(2)}$ be the quotient map.  Let $p_1,\ldots, p_m$ be general
  points in $D^{(2)}$ and suppose for each $i$ that
  $\sigma_D^{-1}(p_i)=\{q_i^1,q_i^2\}$.  Then for any nef class $L\in
  N^1(D^{(2)})_{\mathbb R}$,
$$\epsilon(p_1,\ldots,p_m; D^{(2)},L) = \epsilon( q_1^1,q_1^2,\ldots,q_m^1,q_m^2; D^{\times 2}, \sigma_D^* L).$$
\end{lem}
\begin{proof}
  Let $p_X\colon X \to D^{(2)}$ be the blowup of $D^{(2)}$ at $p_1,\ldots,p_m$ with exceptional divisor $E$, and $p_Y\colon Y\to D^{\times 2}$ be the blowup of $D^{\times 2}$ at $q_1^1,q_2^2,\ldots,q_m^1,q_m^2$ with exceptional divisor $F$.   Then under the induced map $\tilde{\sigma}_D\colon Y\to X$ that lifts $\sigma_D$ we have for any $c>0$ that $\tilde{\sigma}_D^* (p_X^* L - cE) = p_Y^* L-cF$.   Since $\tilde{\sigma}_D$ is surjective this implies $p_X^* L - cE$ is nef if and only if $p_Y^* L-cF$ is nef, which proves the lemma.
\end{proof}




\section{Proofs}\label{sec:proofs}

Let $\CC\to T$ be a family of smooth curves of genus $g$ over a disc
$T$ that develops a node.  By this we mean that the family is proper
and flat, the fibre $C_t$ over $t\in T$ is a smooth curve for $t\neq
0$ and that the fibre $C_0$ over $0\in T$ is an irreducible curve that
has a single node.  We assume further that $\CC$ has a smooth total
space.  Taking the relative second symmetric product of $\CC$ gives a
degeneration of $C^{(2)}$ which when suitably blown up is essentially
the Franchetta degeneration used in
\cite{ciliberto-kouvidakis(99):symmet_produc_a_curve_with_gener_modul}.
Instead of using this we find it easier to consider the fibred
product $\CC\times_T\CC$ (and thus a degeneration of $C\times C$)
and only deal with divisors that are invariant under permuting the
factors.

To this end suppose $p$ is the node in $C_0$ and let $\Y\to \CC
\times_T \CC$ be the blowup of $\CC\times_T\CC$ at $(p,p)$ with
exceptional divisor $E$.  One can easily check by working in local
analytic coordinates that $\Y$ has smooth total space.  Clearly the
fibre of $\Y$ over $t\neq 0$ is $\Y_t=C_t\times C_t$ and the next two
lemmas describe the central fibre of $\Y$.  Denote the normalisation
of $C_0$ by $D$ and let $q,r\in D$ be the preimage of the node $p\in
C_0$.


\begin{lem}\label{lem:descriptionofcentralfibreofY}\ 
  \begin{enumerate}
  \item The central fibre $\Y_0$ has two irreducible components namely
    the exceptional divisor $ E$ which is isomorphic to $\PP^1\times
    \PP^1$ and another we denote by $F$.
  \item The normalisation $\tilde{F}$ of $F$ is the blowup $\pi\colon
    \tilde{F}\to D\times D$ at the four points $(q,q),(q,r),(r,q),(r,r)$
making the natural diagram 
$$
  \begin{CD}
\tilde{F} @>>> F \\
@V\pi VV @VVV\\
D\times D @>>> C_0\times C_0
  \end{CD}  
  $$
  commute.  We denote the exceptional curve in $\tilde{F}$ that sits
  over the point $(s,t)$ by $\tilde{e}_{st}$ and the corresponding
  curve in $F$ by $e_{st}$.
  \item The two components of $\Y_0$ are glued along the four rational
    curves $\{e_{qq},e_{rr}\}$ and $\{e_{qr},e_{rq}\}$ in $F$.  Each
    set consists of a pair of lines in one of the two rulings of $E$.
  \end{enumerate}
\end{lem}
\begin{proof}
  In local analytic coordinates around the node $p\in C_0$ the family
  $\CC$ has the form $xy=t$ in $\mathbb C^2\times T$ where $t$ is the
  parameter on $T$.  Hence locally $\CC\times_T \CC$ is given by
  $x_1y_1=x_2y_2=t$ in $\mathbb C^4\times T$.  Thus a local model for
  $\Y$ is given by the proper transform in the blowup $B\to \mathbb
  C^4\times T$ at the origin.  On the exceptional $\PP^4$ in $B$ we
  pick coordinates $\lambda_1,\lambda_2,\mu_1,\mu_2,\sigma$ such that
  for $t\neq 0$
$$\frac{\lambda_i}{x_i} = \frac{\mu_j}{y_j} = \frac{\sigma}{t} \quad i,j=1,2.$$
Then $E$ is the intersection of $\Y$ with $\PP^4$ and is given by
$$\lambda_1\mu_1=\lambda_2\mu_2 \text{ and } \sigma=0$$
which is a quadric hypersurface in $\PP^3$, and thus $E\simeq
\PP^1\times \PP^1$ as claimed.

Now let $U$ and $V$ be the two components of the normalisation of
$xy=t$ corresponding to $x=0$ and $y=0$ respectively.  Then locally the other
components of $\Y$ are the proper transform of $U\times U,U\times V,V\times U$ and $V\times
V$ given by $x_1=x_2=0$, $x_1=y_2=0$, $y_1=x_2=0$ and $y_1=y_2=0$
respectively.  These are glued along normal crossing curves and the
normalisation $\tilde{F}$ is obtained by pulling them apart.  Thus
$\tilde{F}$ is the blowup of $D\times D$ in the four points as
claimed.

Now the proper transform of $U\times U$ is the blowup at the point
$(p,p)$ and meets $E$ in the line given by $\lambda_1=\lambda_2=0$.
Similarly $V\times V$ meets $E$ in the line $\mu_1=\mu_2$ which is
easily seen to be in the same ruling.  A completely analogous analysis
applies to $U\times V$ and $V\times U$.
\end{proof}

For simplicity denote the numerical class of the curves
$\tilde{e}_{st}$ by the same letter.  Then $$N^1(\tilde{F}) = \pi^*
N^1(D\times D) \oplus \mathbb
Z[\tilde{e}_{qq},\tilde{e}_{rr},\tilde{e}_{qr},\tilde{e}_{rq}].$$
Moreover $N^1(E)$ is a free group of rank $2$ with two generators
$\alpha$ and $\beta$.  We declare that $\alpha$ is the class of the
curve $e_{qq}$ (equivalently of $e_{rr}$) inside $E$ and $\beta$ is
the class of $e_{qr}$ (equivalently of $e_{rq}$).


Consider now the proper transform $\D\subset \Y$ of the diagonal in
$\CC\times_T \CC$.

\begin{lem}\label{lem:diagonalinY}
  The restriction of $\D$ to $E$ has class $\alpha\in N^1(E)$.  The
  pullback of $\D|_F$ to $\tilde{F}$ is the proper transform of
  the diagonal $\Delta_D \subset D\times D$ and thus has class $\pi^*
  \Delta_D - \tilde{e}_{qq} - \tilde{e}_{rr} \in N^1(\tilde{F})$.
\end{lem}
\begin{proof}
  We continue to use the local coordinates introduced in the proof of
  Lemma \ref{lem:descriptionofcentralfibreofY}.  For $t\neq 0$ the
  diagonal is given by $x_1=y_1$ and $x_2=y_2$ and thus meets $E$ in
  the line $\lambda_1=\mu_1$ and $\lambda_2=\mu_2$.  As is easily
  checked this line has class $\alpha$. 

  Now clearly the proper transform of $\D|_{U\times U}$ is, locally,  the
  proper transform of the diagonal $y_1=y_2$ in $U\times U$ and
  similarly for $V\times V$.  Moreover for $\Delta_D$ is disjoint from
  $U\times V$ and $V\times U$ for $t\neq 0$ and thus the pullback of
  $\D|_{F}$ to $\tilde{F}$ has the numerical class as claimed.
\end{proof}

With these preliminaries we are ready to give the proofs of the
theorems stated in the introduction.

\begin{proof}[Proof of Theorem \ref{thm:ampleandseshadristronger}]
  By hypothesis there is a smooth curve $D$ of genus $g-1$ and
  distinct points $q,r\in D$ such that
\begin{equation}
  \epsilon(q+r; D^{(2)}, (a+b)x-b(\delta/2))\ge b\label{eq:hypothesis}.
\end{equation}
By gluing $q$ and $r$ we get a curve $C_0$ with a single node $p$
whose arithmetic genus is $g$ and whose normalisation is $D$.  Let
$\CC\to T$ be a proper flat family of curves of genus $g$ which has a
smooth total space, smooth general fibre and whose central fibre is
$C_0$.  Let $\Y\to \CC\times_T\CC$ be the blowup at $(p,p)$.  We will
use the notation introduced in Lemmas
\ref{lem:descriptionofcentralfibreofY} and \ref{lem:diagonalinY} so
the central fibre $\Y_0$ has two components $E$ and $F$, and the
normalisation $\tilde{F}$ of $F$ is the blowup $\pi\colon \tilde{F}\to
D\times D$ at the four points $(q,q),(q,r),(r,q),(r,r)$.

Fix a line bundle $L$ on $\CC/\CC_0$ that has degree 1 on each of the
fibres.  As $\CC$ is assumed to be smooth and $\CC_0$ is irreducible,
$L$ extends uniquely to a line bundle $L'$ on all of $\CC$.  Take a
meromorphic section of $L'$ whose support does not contain the node
$p$ of $C_0$.  By shrinking $T$ if necessary we may assume furthermore
that the support of $s$ does not contain any fibre $\CC_t$, and we
write this support as $\sum_i a_i D_i$ for some divisors $D_i\subset \CC$ that do
not contain $p$.  For each $i$ define a divisor on $\Y$ by
$$\G_i = \{ (u, v)\in \Y_t \text{ for some } t \text{ and either } u\in D_i \text{ or } v\in D_i\},$$
and set
$$\G = \sum_i a_i \G_i.$$

Clearly $\G$ is invariant under the natural involution on $\Y$.  In
fact for for all $t\neq 0$ the numerical class of $\G|_{\Y_t}$ is the
pullback of $x$ under $\sigma_{C_t} \colon C_t\times C_t\to
C_t^{(2)}$.  Moreover $\G$ is trivial along $E$ and and the numerical
class of the pullback of $\G|_F$ to $\tilde{F}$ has class $\pi^* \sigma_D^*
x$.

Let $\D\subset \Y$ be the proper transform of the diagonal in
$\CC\times_T \CC$ as in Lemma \ref{lem:diagonalinY} and define an
$\mathbb R$-divisor on $\Y$ by
\begin{equation}
 \H = (a + b) \G - b (\D + E).\label{eq:definitionofthedivisorH}
\end{equation}
We claim that $\H|_{\Y_0}$ is nef.  To see this note first from
(\ref{lem:descriptionofcentralfibreofY}, \ref{lem:diagonalinY}) that
$E|_E$ has class $-(2\alpha+2\beta)$ and $\D|_E$ has class $\alpha$.
Thus $\H|_E$ has class $b (\alpha+2\beta)$ which is clearly nef as
both $\alpha$ and $\beta$ are nef and $b>0$.  To show that $\H|_F$ is
nef consider its pullback to $\tilde{F}$ which using
(\ref{lem:descriptionofcentralfibreofY},\ref{lem:diagonalinY}) has
numerical class
\begin{multline}\label{eq:classofHrestrictedtoFtilde}
    (a+b) \pi^*\sigma_D^* x - b(\pi^* \Delta_D - \tilde{e}_{qq} - \tilde{e}_{rr}) - b(\tilde{e}_{qq}+\tilde{e}_{rr}+\tilde{e}_{rq}+\tilde{e}_{rq}) \\
    =\pi^* \sigma_D^* \left((a+b)x - b(\delta/2)\right) -
    b(\tilde{e}_{rq} + \tilde{e}_{rq})
\end{multline}
since $\sigma_D^*(\delta/2)=\Delta_D$.  Now $\sigma_D^{-1}(q+r) =
\{(q,r),(r,q)\}$ so \eqref{eq:hypothesis} and Lemma \ref{lem:seshadriquotient} imply
$$\epsilon( (q,r),(r,q); D\times D, \sigma_D^* ((a+b)x - b(\delta/2)))\ge b.$$ 
But by \eqref{eq:classofHrestrictedtoFtilde} this means exactly that
the pullback of $\H|_{F}$ to $\tilde{F}$ is nef and thus $\H|_F$ is
nef as well.

Hence $\H|_{\Y_0}$ is nef and by semicontinuity so is $\H|_{\Y_t}$
for very general $t$.  But $\H_{\Y_t}$ has class
$\sigma_{C_t}^*((a+b)x-b(\delta/2))$ so $(a+b)x-(\delta/2)\in
N^1(C_t^{(2)})$ is nef for very general $t$ which proves that
$\tau(C_t)\le \frac{a}{b}$.  By \eqref{lem:jumpsdown} the same
inequality holds for any very general curve of genus $g$.
\end{proof}

\begin{proof}[Proof of Theorem \ref{thm:seshadrimulti}]
  The fact that $(a+b)x-b(\delta/2)\in N^1(D^{(2)})$ is nef comes from
  Theorem \ref{thm:ampleandseshadristronger} since
  $\epsilon_{r}(\cdot) \le \epsilon_1(\cdot)$.  Essentially the result
  we want comes from the degeneration $\Y$ described in the proof of
  Theorem \ref{thm:ampleandseshadristronger} and semicontinuity of
  Seshadri constants in families.

  To describe this more explicitly we continue the notation from the
  above proof.  Pick $m$ sections $s_1,\ldots, s_m$ of $\CC\times_T
  \CC\to T$ that meet $C_0\times C_0$ at $m$ very general points (so
  in particular these points are not equal to $(p,p)$).  Let $V_i$ be
  the image of $s_i$ and $V_i'$ be the image of $V_i$ under the
  involution.  We denote the proper transform of $V_i$ and $V_i'$ in
  $\Y$ by $W_i$ and $W_i'$ and let $W=\bigcup_i W_i \cup W_i'$.  By
  shrinking $T$ if necessary we may assume that $W$ meets each fibre
  $\Y_t=C_t\times C_t$ at a collection of $2m$ distinct points.  Note
  that in the central fibre $\Y_0$ these points are all in the
  component $F$.

  Now let $\pi\colon \Y'\to \Y$ be the blowup along $W$ with
  exceptional divisor $E'$.  The central fibre of $\Y'$ has components
  $E$ and $F'$ where the normalisation $\tilde{F'}$ of $F'$ is the
  blowup of $\tilde{F}$ at $m$ very general points and their image
  under the involution.  Thus $\tilde{F'}$ is the blowup of $D\times
  D$ at the points $(q,q),(q,r),(r,q),(r,r)$ and at a further $2m$
  points.  

  Set $c=\epsilon_m(D^{(2)}, (a+b)x-(\delta/2))$ and consider the
  $\mathbb R$-divisor
$$ H' = \pi^* H - c E'$$
where $H$ is the divisor defined in
\eqref{eq:definitionofthedivisorH}.  Then exactly as in the proof of
Theorem \ref{thm:ampleandseshadristronger}, $H'|_E$ is nef and the
hypothesis on the Seshadri constant and Lemma
\ref{lem:seshadriquotient} imply that the pullback of $H'|_{F'}$ to
$\tilde{F'}$ is also nef.  Thus $H'|_{\Y'_t}$ is nef for very general
$t$ and using \eqref{lem:seshadriquotient} once again proves the
theorem.
\end{proof}

\begin{proof}[Proof of Corollary
  \ref{cor:ampleintermsofseshadriofpp2}]
  Fix $g\ge 1$.  Let $h$ be the class of the hyperplane in $\PP^2$ and
  set $b=\epsilon_g(\PP^2,h)$.  The second symmetric product of a
  genus 0 curve is $\PP^2$ and under this identification
  $x=(\delta/2)=h$.  We have
   $$\epsilon_{g}((\PP^1)^{(2)},(1+b)x-b(\delta/2)) = \epsilon_g(\PP^2,h)=b.$$ Thus repeated
   use of Theorem \ref{thm:seshadrimulti} yields for a very general
   curve $D$ of genus $g-1$
$$\epsilon_1(C^{(2)},(1+b)x-b(\delta/2))\ge b,$$
and so the result follows from Theorem
\ref{thm:ampleandseshadristronger}.
\end{proof}

\section{Application to the case $g=5$}\label{sec:applicationgenus5}

We now prove that if $C$ is a very general curve of genus $5$ then
$\tau(C)\le 16/7$.  This is done by estimating the Seshadri constant
at a very general point $p$ of $D^{(2)}$ where $D$ is a very general
curve of genus $4$.  By \eqref{sec:lowgenus} we know that $\tau(D)=2$.
Set $a=16,b=7$ and $L=(a+b)x-b(\delta/2)\in N^1(D^{(2)})$ which is
ample.  By
\eqref{thm:ampleandseshadristronger} it is sufficient to show the following \\

\noindent {\it Claim:  If $p$ is a very general point in $D^{(2)}$ then 
  \begin{equation}
\epsilon(p; D^{(2)}, L)\ge b=7.\label{eq:claimseshadrig4}
\end{equation}
}

The proof of the claim will use the ideas of Ein-Lazarsfeld
\cite{ein-lazarsfeld:93:seshad_const_smoot_surfac}.  Rather than using
the main result of that paper we get an improvement by using their
techniques and special properties of the particular surface $D^{(2)}$.
In particular we will need the following lemma.

\begin{lem}[Ein-Lazarsfeld  \cite{ein-lazarsfeld:93:seshad_const_smoot_surfac}]\label{lem:einlazarsfeld}
Let $X$ be a smooth surface and $L$ be an integral ample line bundle
(or class) on $X$.  Suppose $\{p_t\in E_t\}_{t\in T}$ is a
one-parameter family consisting of a point $p_t$ in a curve $E_t\subset
X$ such that $\mult_{p_t} E_t\ge m$ for all $t$.  Suppose in addition
that $E=E_0$ is reduced and irreducible and moreover that the
Kodaira-Spencer class of this family is non-zero.  Then $E^2\ge
m(m-1)$.
\end{lem}

Now suppose $E\subset D^{(2)}$ is a reduced irreducible curve passing
through a very general point $p$ with
numerical class $(n+\gamma) x - \gamma(\delta/2)$.

\noindent {\it Claim: We have $L.E\ge 7$.  Moreover if $(n,\gamma)
  \notin \{(1,0),(3,1),(5,2)\}$
  then $L.E\ge 4b=28$.  }\\

To see this note that $L.E = an - 4b\gamma=16n-28\gamma$ so certainly
if $\gamma\le 0$ then $L.E\ge 28$ unless $(n,\gamma)=(1,0)$.  So
suppose that $\gamma>0$.  As $\tau(D)=2$ we must have $n\ge 2\gamma$
with equality if and only if $E$ has zero self intersection.  

Now for each fixed $\gamma\ge 0$ there are at most finitely many
irreducible curves $E$ with numerical class
$(n+\gamma)x-\gamma(\delta/2)$ and self-intersection zero.  Since $p$
is assumed to be very general we may assume there is no irreducible
curve of zero self-intersection through $p$, and so we in fact have
$n\ge 2\gamma+1$.  Then it is easily checked that $L.E\ge 7$ and $L.E\ge
28$ except when
$(n,\gamma)\in \{(3,1),(5,2)\}$.  \\

\noindent We now finish the proof of \eqref{eq:claimseshadrig4}.
Suppose for contradiction that $\epsilon(p; D^{(2)},L)<7$ for a very
general $p\in D^{(2)}$.  Then through a very general point $p$ there
exists a reduced irreducible curve $E$ with $m=\mult_p E$ and
\begin{equation}
\frac{L.E}{m} < b=7.\label{eq:seshadriexceptional}
\end{equation}
As in \cite{ein-lazarsfeld:93:seshad_const_smoot_surfac} the
collection of pairs $(p,E)$ consisting of a point $p$ in an
irreducible curve $E$ such that $\mult_p(E)>L.E/7$ consists of a
countable collection of algebraic families and the proof will be
completed by showing that any such family with $p$ a very general
point is discrete.

To this end suppose for contradiction that there is a family $\{p_t\in
E_t\}_{t\in T}$ with $E_t$ reduced and irreducible and $\mult_{p_t}
E_t > L.E_t/7$ for all $t$. Set $(p,E)=(p_0,E_0)$.  Since $p$ is
assumed to be very general we have from the above that $L.E_0\ge 7$ so
$\mult_{p_t}>1$.  Since $\mult_y E_t=1$ for a general point $y$ of
$E_t$ we deduce that the curves $E_t$ are moving in a non-trivial
family and thus from \eqref{lem:einlazarsfeld} we have
that $E^2\ge m(m-1)$. \\

\noindent{\it Case 1: $(n,\gamma) = (1,0)$.} If $E$ is not the irreducible curve $x_p=\{p+q | q\in C\}$ then $m\le E.x=1$ so $m=1$ which is impossible by \eqref{eq:seshadriexceptional}. On the other hand if $E=x_p$ then $m=1$ and $L.E=a$ which again is absurd. \\

\noindent{\it Case 2: $(n,\gamma) = (3,1)$ (resp.\ $(n,\gamma) =
  (5,2)$)}.  As $m(m-1)\le E^2 =5$ (resp.\ $m(m-1)\le E^2=9$) we have $m\le 2$ (resp.\ $m\le 3$).
But this implies $\frac{L.E}{m} \ge \frac{3a-4b}{2}=10\ge b$ (resp.\ $\frac{L.E}{m}\ge
\frac{5a-8b}{3}=8\ge b$) which in both cases is impossible by \eqref{eq:seshadriexceptional}.\\

\noindent{\it Case 3: $(n,\gamma)\neq (1,0),(3,1),(5,2)$}.  Here we
follow \cite{ein-lazarsfeld:93:seshad_const_smoot_surfac} but first
note that by the previous claim in this case $L.E\ge 4b$ which by
\eqref{eq:seshadriexceptional} implies that $m\ge 5$.  Again from
\eqref{eq:seshadriexceptional} we have $L.E<7m$ so $L.E\le 7m-1$.
Thus by the Hodge index theorem
$$ m(m-1) \le E^2 \le \frac{(L.E)^2}{L^2} \le \frac{(7m-1)^2}{60}$$
which is impossible for $m\ge 5$.

\bibliography{biblio} 

\begin{thebibliography}{10}

\bibitem{arbarello_cornalba_griffiths(85):geomet_algeb_curves}
E.~Arbarello, M.~Cornalba, P.~A. Griffiths, and J.~Harris.
\newblock {\em Geometry of algebraic curves. {V}ol. {I}}, volume 267 of {\em
  Grundlehren der Mathematischen Wissenschaften [Fundamental Principles of
  Mathematical Sciences]}.
\newblock Springer-Verlag, New York, 1985.

\bibitem{bauer(99):seshad_const_algeb_surfac}
T.~Bauer.
\newblock Seshadri constants on algebraic surfaces.
\newblock {\em Math. Ann.}, 313(3):547--583, 1999.

\bibitem{campana_peternell(90):algeb_ample_cone_projec_variet}
F.~Campana and T.~Peternell.
\newblock Algebraicity of the ample cone of projective varieties.
\newblock {\em J. Reine Angew. Math.}, 407:160--166, 1990.

\bibitem{ciliberto-kouvidakis(99):symmet_produc_a_curve_with_gener_modul}
C.~Ciliberto and A.~Kouvidakis.
\newblock On the symmetric product of a curve with general moduli.
\newblock {\em Geom. Dedicata}, 78(3):327--343, 1999.

\bibitem{ein_kuechle_lazarsfeld(95):local_posit_ample_line_bundl}
L.~Ein, O.~K{\"u}chle, and R.~Lazarsfeld.
\newblock Local positivity of ample line bundles.
\newblock {\em J. Differential Geom.}, 42(2):193--219, 1995.

\bibitem{ein-lazarsfeld:93:seshad_const_smoot_surfac}
L.~Ein and R.~Lazarsfeld.
\newblock Seshadri constants on smooth surfaces.
\newblock {\em Ast\'erisque}, 218:177--186, 1993.
\newblock Journ\'ees de G\'eom\'etrie Alg\'ebrique (Orsay, 1992).

\bibitem{harbourne_roe:05:multip_seshad_const_bbb_p_sp_2}
B.~Harbourne and J.~Ro{\'e}.
\newblock Multipoint {S}eshadri constants on {$\Bbb P\sp 2$}.
\newblock {\em Rend. Sem. Mat. Univ. Politec. Torino}, 63(1):99--102, 2005.

\bibitem{kouvidakis:93:divis_symmet_produc_curves}
A.~Kouvidakis.
\newblock Divisors on symmetric products of curves.
\newblock {\em Trans. Amer. Math. Soc.}, 337(1):117--128, 1993.

\bibitem{kouvidakis:02:some_resul_morit_their_applic}
A.~Kouvidakis.
\newblock On some results of {M}orita and their application to questions of
  ampleness.
\newblock {\em Math. Z.}, 241(1):17--33, 2002.

\bibitem{lazarsfeld(04):posit_in_algeb_geomet}
R.~Lazarsfeld.
\newblock {\em Positivity in algebraic geometry. {I}}, volume~48 of {\em
  Ergebnisse der Mathematik und ihrer Grenzgebiete. 3. Folge.}
\newblock Springer-Verlag, Berlin, 2004.
\newblock Classical setting: line bundles and linear series.

\bibitem{oguiso(02):seshad_const_in_a_famil_surfac}
K.~Oguiso.
\newblock Seshadri constants in a family of surfaces.
\newblock {\em Math. Ann.}, 323(4):625--631, 2002.

\bibitem{pacienza_polizzi:01:a_degen_symmet_produc_a}
G.~Pacienza and F.~Polizzi.
\newblock On a degeneration of the symmetric product of a curve with general
  moduli.
\newblock {\em Matematiche (Catania)}, 56(2):297--307 (2003), 2001.
\newblock PRAGMATIC, 2001 (Catania).

\bibitem{strycharz_szemberg_szemberg(04):remar_nagat_conjec}
B.~Strycharz-Szemberg and T.~Szemberg.
\newblock Remarks on the {N}agata conjecture.
\newblock {\em Serdica Math. J.}, 30(2-3):405--430, 2004.

\bibitem{szemberg_tutaj_gasinska(02):gener_blow_ups_projec_plane}
T.~Szemberg and H.~Tutaj-Gasi{\'n}ska.
\newblock General blow-ups of the projective plane.
\newblock {\em Proc. Amer. Math. Soc.}, 130(9):2515--2524 (electronic), 2002.

\bibitem{xu(94):curves_in_f_p_sp}
G.~Xu.
\newblock Curves in {${\bf P}\sp 2$} and symplectic packings.
\newblock {\em Math. Ann.}, 299(4):609--613, 1994.

\end{thebibliography}
\vspace{4mm}

\noindent Julius Ross, Department of Mathematics,\\
 Columbia University, New York, NY 10027. USA. \\

\end{document}